\documentclass[11pt]{article}

\usepackage[a4paper]{anysize}\marginsize{2.5cm}{2.5cm}{1.3cm}{2cm}
\pdfpagewidth=\paperwidth \pdfpageheight=\paperheight
\usepackage{amsfonts,amssymb,amsthm,amsmath,eucal}
\usepackage{pgf}
\usepackage{bbm}
\usepackage{tikz} 
\usepackage{mathrsfs}
\usetikzlibrary{arrows}
\usepackage{color}
\usepackage{float}
\usepackage{subfigure}
\usepackage{caption}
\usepackage{graphicx}
\usepackage{url}
\usepackage[most]{tcolorbox}
\usepackage{hyperref}
\usepackage{cleveref}

\definecolor{qqqqff}{rgb}{0.,0.,0.}

\theoremstyle{plain}
\newtheorem{thm}{Theorem}[section]
\newtheorem{theorem}[thm]{Theorem}

\newtheorem{lemma}[thm]{Lemma}
\newtheorem{proposition}[thm]{Proposition}
\newtheorem{corollary}[thm]{Corollary}
\newtheorem{conjecture}[thm]{Conjecture}

\theoremstyle{definition}
\newtheorem{definition}[thm]{Definition}
\newtheorem{remark}[thm]{Remark}

\newtheorem{notation}[thm]{Notation}

\newtheorem{thevarthm}[thm]{\varthmname}

\newenvironment{varthm*}[1]{\trivlist\item[]{\bf #1.}\it}{\endtrivlist}

\renewcommand\geq{\geqslant}

\renewcommand\leq{\leqslant}

\newcommand\be{\begin{eqnarray*}}
\newcommand\ee{\end{eqnarray*}}

\newcommand\Z{\mathbb Z}

\newcommand\PP{\mathbb P}

\newcommand\calo{{\mathcal O}}
\newcommand\cali{{\mathcal I}}

\newcommand\newop[2]{\def#1{\mathop{\rm #2}\nolimits}}
\newop\log{log}
\newop\ord{ord}
\newop\Gal{Gal}
\newop\SL{SL}
\newop\Bl{Bl}
\newop\mult{mult}
\newop\mass{mass}
\newop\div{div}
\newop\codim{codim}
\newop\sing{sing}
\newop\Zeroes{Zeroes}

\newcommand{\jd}{\frac{1}{2}}

\newtheorem{innercustomgeneric}{\customgenericname}
\providecommand{\customgenericname}{}
\newcommand{\newcustomtheorem}[2]{%
  \newenvironment{#1}[1]
  {%
   \renewcommand\customgenericname{#2}%
   \renewcommand\theinnercustomgeneric{##1}%
   \innercustomgeneric
  }
  {\endinnercustomgeneric}
}

\newcustomtheorem{customscheme}{Scheme}
\newcustomtheorem{customlemma}{Lemma}

\makeatletter
\newcommand\setboxcounter[2]{\setcounter{tcb@cnt@#1}{#2}}
\makeatother

\newtcolorbox[auto counter,number within=section]{schemeM}[2][]{
  enhanced, hbox, 
  colback=blue!5!white,
  colframe=blue!75!black,
  fonttitle=\scshape,
  title={Scheme~\thetcbcounter. #2 #1},
  #1
}

\newtcolorbox[auto counter,number within=section]{scheme}[1][]{
  enhanced, hbox, 
  colback=blue!5!white,
  colframe=blue!75!black,
  fonttitle=\scshape,
  title={Scheme \thetcbcounter},
  #1
}

\newtcolorbox[]{schemeN}[1][]{
  enhanced, hbox, 
  colback=blue!5!white,
  colframe=blue!75!black,
  fonttitle=\scshape,
  title={Scheme \thetcbcounter},
  #1
}

\newcommand{\initial}[1]{
\begin{scheme}
#1
\end{scheme}
}

\newcommand{\initialN}[1]{
\begin{schemeN}
#1
\end{schemeN}
}

\newcommand{\reduction}[3]{\begin{tcolorbox}[colback=red!5!white,colframe=red!75!black,height=1cm,hbox]
  {#1}
\end{tcolorbox}
\vspace{.1cm}

\begin{tcolorbox}[colback=white!5!white,colframe=white!75!white,width=.15\textwidth,height=1cm,before=\hfill,after=\hfill]
  $\big\downarrow$ $\longrightarrow$
\end{tcolorbox}
\begin{tcolorbox}[colback=green!5!white,colframe=green!75!black,width=.75\textwidth,before=\hfill,height=1cm,after=\hfill]
  \centering{#2}
\end{tcolorbox}

\begin{scheme}
#3
\end{scheme}}

\newcommand{\reductionL}[4]{\begin{tcolorbox}[colback=red!5!white,colframe=red!75!black,height=1cm,hbox]
  {#1}
\end{tcolorbox}
\vspace{.1cm}

\begin{tcolorbox}[colback=white!5!white,colframe=white!75!white,width=.15\textwidth,height=1cm,before=\hfill,after=\hfill]
  $\big\downarrow$ $\longrightarrow$
\end{tcolorbox}
\begin{tcolorbox}[colback=green!5!white,colframe=green!75!black,width=.75\textwidth,before=\hfill,height=1cm,after=\hfill]
  \centering{#2}
\end{tcolorbox}

\begin{scheme}[label=#4]
#3
\end{scheme}}
\newcommand{\reductionN}[3]{\begin{tcolorbox}[colback=red!5!white,colframe=red!75!black,height=1cm,hbox]
  {#1}
\end{tcolorbox}
\vspace{.1cm}

\begin{tcolorbox}[colback=white!5!white,colframe=white!75!white,width=.15\textwidth,height=1cm,before=\hfill,after=\hfill]
  $\big\downarrow$ $\longrightarrow$
\end{tcolorbox}
\begin{tcolorbox}[colback=green!5!white,colframe=green!75!black,width=.75\textwidth,before=\hfill,height=1cm,after=\hfill]
  \centering{#2}
\end{tcolorbox}

\begin{schemeN}
#3
\end{schemeN}}

\newcommand{\reductionM}[4]{\begin{tcolorbox}[colback=red!5!white,colframe=red!75!black,height=1cm,hbox]
  {#1}
\end{tcolorbox}
\vspace{.1cm}

\begin{tcolorbox}[colback=white!5!white,colframe=white!75!white,width=.15\textwidth,height=1cm,before=\hfill,after=\hfill]
  $\big\downarrow$ $\longrightarrow$
\end{tcolorbox}
\begin{tcolorbox}[colback=green!5!white,colframe=green!75!black,width=.75\textwidth,before=\hfill,height=1cm,after=\hfill]
  \centering{#2}
\end{tcolorbox}

\begin{schemeM}{#4}
#3
\end{schemeM}}

\def\keywordname{{\bfseries Keywords}}%
\def\keywords#1{\par\addvspace\medskipamount{\rightskip=0pt plus1cm
\def\and{\ifhmode\unskip\nobreak\fi\ $\cdot$
}\noindent\keywordname\enspace\ignorespaces#1\par}}
\def\subclassname{{\bfseries Mathematics Subject Classification
(2000)}\enspace}
\def\subclass#1{\par\addvspace\medskipamount{\rightskip=0pt plus1cm
\def\and{\ifhmode\unskip\nobreak\fi\ $\cdot$
}\noindent\subclassname\ignorespaces#1\par}}

\begin{document}
\vspace{-2cm}

\title{Postulation of lines in $\PP^3$ revisited}
\author{Marcin Dumnicki, Mikołaj Le Van, Grzegorz Malara, Tomasz Szemberg,\\ Justyna Szpond and Halszka Tutaj-Gasi\'nska}
\date{}
\maketitle
\thispagestyle{empty}
\parskip=.05cm

\begin{abstract}
The purpose of the present note is to provide a new proof ot the well-known result due to Hartshorne and Hirschowitz to the effect that general lines in projective spaces have \emph{good postulation}. Our approach uses specialization to a hyperplane and thus opens door to study postulation of general codimension $2$ linear subspaces in projective spaces.
\subclass{14C20; 14N05; 14N15; 14D06; 13H15; 13D40}
\end{abstract}

\section{Introduction}
One of classical problems in algebraic geometry is to study the dimension of linear systems with imposed base loci. In this note we are interested in the situation that the base locus is a finite union of general codimension $2$ linear subspaces in a projective space $\PP^N$. Linear systems of this kind appear in various guises across algebraic geometry, commutative algebra and combinatorics, see e.g. \cite{Ida_Crelle_1990}, \cite{Rice23}, \cite{Advances2014}, \cite{Derksen_Sidman_2002}, \cite{Schenck_Sidman_2013}.  

Let us recall the following general notion of \emph{good postulation}.
\begin{definition}[Good postulation]
A closed subscheme $Y$ in a projective space $\PP^N$ is said to have \emph{good postulation}, if the restriction map
\begin{equation}\label{eq: restriction}
    H^0(\PP^N,\calo_{\PP^N}(d))\to H^0(Y,\calo_Y(d))
\end{equation}
has maximal rank for all $d\geq 1$.
\end{definition}
Expressed in a slightly different language \emph{good postulation} is equivalent to saying that the subscheme $Y$ imposes \emph{the expected number of conditions} on forms of arbitrary degree. 

Our research is driven by the following conjecture.
\begin{conjecture}[Codimension 2 good postulation]\label{conj: good postulation}
A finite union of general codimension $2$ linear subspaces in $\PP^N$ has good postulation.    
\end{conjecture}
This Conjecture should be viewed as going in a complimentary direction when compared with the predictions of Carlini, Catalisano and Geramita in \cite{CCG2010} that a finite union of \emph{disjoint} linear subspaces imposes independent conditions on forms of arbitrary degree.

Conjecture \ref{conj: good postulation} is obviously true in $\PP^2$. It has been proved in $\PP^3$ by Hartshorne and Hirschowitz in \cite{HarHir82}. Their strategy was to specialize some of the lines in a ruling of a smooth quadric in $\PP^3$. There are no hypersurfaces with analogous properties, i.e., such that they contain any number of codimension $2$ subspaces not inflicting extra intersections among them, in higher dimensional projective spaces. This was a considered a serious obstacle to generalize the result of Hartshorne and Hirschowitz to higher dimensions, see e.g. \cite{Ballico_Few_Lines2011}, \cite{Tovar_Planes_P5_1985} for some initial attempts in that direction.

Aladpoosh and Catalisano proposed in \cite{AlaCat21} a new proof of the good postulation of lines in $\PP^3$, which does not use degeneration to a smooth quadric. Instead they introduce degenerations to non-reduced subschemes of a rather complicated structure. In contrast, our approach seems more natural as we degenerate some lines to a plane $\PP^2$ in $\PP^3$. We hope that the method proposed here can be easier adapted to higher dimensional spaces, see \cite{Chmiel21} for some preparatory facts concerning planes in $\PP^4$.

We provide in \cite{SingularFiles} software, which not only allows to check all calculations carried out below, but examines all reduction pathes allowing sometimes for easier but sporadic arguments. The readers are encourage to experiment with our programs.

\section{The main result and an outline of the proof}\label{sec: main result}
Our main result is the following Theorem, which implies Conjecture \ref{conj: good postulation} as we shall see shortly. The number $\delta_{i,j}$ is the Kronecker symbol equal $1$ if $i=j$ and $0$ otherwise.
\begin{theorem}[Main Theorem]\label{thm: main}
Let $d=3k+\varepsilon$ with $\varepsilon\in\left\{0,1,2\right\}$ and let 
$$Y(\ell,p)$$
be a subscheme of $\PP^3$ consisting of
$\ell=\frac{1}{2}(3k+5+2\varepsilon)k+1+\varepsilon$ general lines and $p=(k+1)\cdot\delta_{2,\varepsilon}$ collinear points (in a line general with respect to the other $\ell$ lines).
Then the map
\begin{equation}\label{eq: main theorem}
    H^0(\PP^N,\calo_{\PP^N}(d))\to H^0(Y(\ell,p),\calo_{Y(\ell,p)}(d))
\end{equation}
is an isomorphism.
\end{theorem}
The numbers in the above Theorem are chosen so that the spaces on both sides have the same dimension. Now, if the map \eqref{eq: main theorem} is an isomorphism, its kernel $H^0(\PP^3,\mathcal{I}_Y(d))$ is trivial, which means that there is no nonzero form of degree
$d$ vanishing on $Y$. Hence, keeping $d$ fixed, for any number $s>\ell$ there is no such form, so that the map is an injection in all these cases. On the other hand, for $s<\ell$ the restriction map must be surjective as otherwise it would fail to be surjective for the scheme in the Main Theorem.

We prove Theorem \ref{thm: main} in the subsequent sections. The idea is to show that subschemes $Y(\ell,p)$ can be specialized to subschemes $B(L,C,(A,B))$ consisting of $L$ general lines, $C$ crosses and $\binom{A}{2}-B$ points forming a $P(A,B)$ sub-star with suitable integers $A,B$, see Notation \ref{not: sub-star}. This is proved in Proposition \ref{prop: initial specialization}.

Then, for schemes of type $B(L,C,(A,B))$ Lemma \ref{lem: reduction} provides induction steps reducing the number of lines and the degree of studied forms lying out the basis for an induction argument.

The initial steps for the induction argument are presented in Section \ref{sec: preinduction cases}.

During the reductions, certain trace subschemes in a projective plane arise. We address their postulation in Section \ref{sec: postulation P2}.

\section{The main tools}\label{sec: main tools}
Let $Y$ be a closed subscheme in $\PP^N$. We denote by $\cali_Y$ its ideal sheaf and by $h_Y$ its Hilbert function, i.e.,
$$h_Y(d)=h^0(\PP^N,\calo_{\PP^N}(d))-h^0(Y,\cali_Y(d)).$$
Our approach builds upon the \emph{semicontinuity theorem} (see \cite[Theorem III.12.8]{Hartshorne}), which asserts that in a flat family of subschemes, the dimension of cohomology groups is upper semicontinuous as a function on the parameter space of the family. This is applied so that demonstrating the Hilbert function pattern for a specific subscheme of a given type suffices to establish it for a general subscheme of that type.

In this work, we are primarily interested in the case where $Y$ is a union of general lines in $\PP^3$. The strategy outlined above is applied by specializing some of the lines to a plane $H \subset \PP^3$. 
We then define the \emph{trace scheme} $Y_{tr}$ as the scheme theoretic intersection $Y\cap H$ and we define the \emph{residual scheme} $Y_{res}$ of $Y$ with respect to $H$ as the scheme given by the ideal sheaf $\mathcal{I}_Y : \mathcal{I}_H$. The following exact sequence of sheaves
\begin{equation}\label{eq: ciag dokladny}
0\to \cali_{Y_{res}}(-H)\to \cali_Y\to \cali_{Y_{tr},H}\to 0    
\end{equation}
yields the fundamental relation between the Hilbert functions of $Y$, $Y_{tr}$ and $Y_{res}$, known as the \emph{Castelnuovo's inequality}:
\begin{equation}\label{eq: Castelnuovo}
h_Y(d)\geq h_{Y_{res}}(d-1)+h_{Y_{tr},H}(d),    
\end{equation}
where $h_{Y_{tr},H}$ is computed for $Y_{tr}\subset\PP^2=H$.
Castelnuovo's inequality allows to apply induction on the number of lines and degree. 
Twisting the exact sequence \eqref{eq: ciag dokladny} with $\calo_{\PP^3}(d)$ we obtain the following convenient form the Castelnuovo's inequality, which is crucial for our reasoning.
\begin{theorem}[Castelnuovo]\label{thm: Castelnuovo}
If the linear system of forms of degree $d-1$ in $\PP^3$ vanishing along $Y_{res}$ is empty (i.e. $H^0(\PP^3,\cali_{Y_{res}}(d-1))=0$) and the linear system of forms of degree $d$ in $H=\PP^2$ vanishing along $Y_{tr}$ is empty (i.e. $H^0(H,\cali_{Y_{tr}}(d))=0$), then there does not exist any nonzero form of degree $d$ in $\PP^3$ vanishing along $Y$ (i.e. $H^0(\PP^3,\cali_{Y}(d))=0$).
\end{theorem}

If we can show that the schemes $Y_{res}$ and $Y_{tr}$ impose independent conditions on forms of degree $d-1$ in $\PP^3$ and degree $d$ in $H=\PP^2$ respectively, then we conclude that $Y$ imposes independent conditions on forms of degree $d$ in $\PP^3$. Sometimes it is not possible to split $Y$ in such a way that the schemes $Y_{res}$ and $Y_{tr}$ impose independent conditions. This is circumvent by allowing some other, less obvious, specializations, which are listed below.

We describe now the residual and trace subschemes for specializations occuring in the sequel. 

\textbf{Line-type specialization} of $k$ general lines to form a star configuration $Z$ and putting $Z$ in $H$. In this case in $H$ we obtain as the trace $Z_{tr}$ the star configuration of $k$ lines (see \cite{StarConfigurations} for details). The residual scheme $Z_{res}$ consists then of $\binom{k}{2}$ coplanar points, the double points of $Z_{tr}$.

\textbf{Sundial-type specialization} of 2 lines to intersect in a sundial $Z$ and putting the singular point of $Z$ in $H$. In this case $Z_{tr}$ is a double point in $H$. The residual scheme $Z_{res}$ is a \emph{cross}, which is an equidimensional curve. See \cite{Sundials} for details how a sundial is obtained as a limit of a flat family of two skew lines.

\textbf{Half-cross-type specialization} of one of the lines in a cross $Z$ to $H$. In this case $Z_{tr}$ is the line in $H$ and $Z_{res}$ is the remaining line.

\textbf{Cross-type specialization} of a cross $Z$ to $H$. In this case $Z_{tr}$ is the cross. The residual scheme $Z_{res}$ is empty.

In the sequel we want to perform various types of specializations at the same type. To this end it is convenient to introduce the following notation.
\begin{notation}
We use the symbol
$$S(\ell,s,h,c)$$
to denote a specialization consisting of 
\begin{itemize}
    \item $\ell$ Line-type specializations;
    \item $s$ Sundial-type specializations;
    \item $h$ Half-cross-type specializations and
    \item $c$ Cross-type specialization.
\end{itemize}
\end{notation}
We now introduce additionally the following notation for convenience.
\begin{notation}[Sub-star set of points]\label{not: sub-star}
We use the symbol $P(a,b)$ to denote a subset of the singular points of a star configuration of $a$ lines with $b$ points omitted in such a way that every configuration line contains at most one omitted point. 
\end{notation}
In particular it is always $b\leq a/2$ and the set $P(a,b)$ is the union of $\binom{a}{2}-b$ points. Figure \ref{fig: a,b} illustrates a sub-star $P(4,2)$.
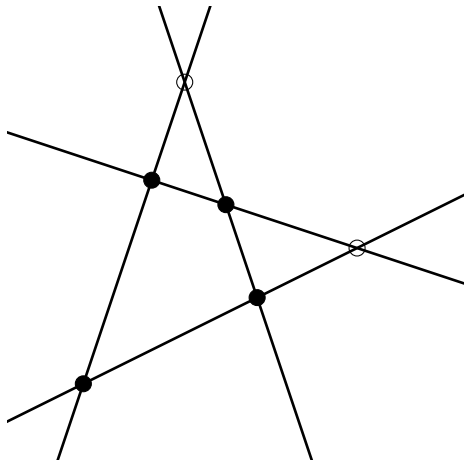
\begin{figure}[h!]
    \centering
\begin{tikzpicture}[line cap=round,line join=round,>=triangle 45,x=1cm,y=1cm]
\clip(-1,-1) rectangle (5,5);
\draw [line width=1pt,domain=-10.373944253023616:15.22960498428316] plot(\x,{(-0--3*\x)/1});
\draw [line width=1pt,domain=-10.373944253023616:15.22960498428316] plot(\x,{(--8-3*\x)/1});
\draw [line width=1pt,domain=-10.373944253023616:15.22960498428316] plot(\x,{(-0--0.5*\x)/1});
\draw [line width=1pt,domain=-10.373944253023616:15.22960498428316] plot(\x,{(--3-0.3333333333333333*\x)/1});
\begin{scriptsize}
\draw [fill=black] (0,0) circle (3pt);
\draw [fill=black] (2.2857142857142856,1.1428571428571428) circle (3pt);
\draw [fill=black] (0.9,2.7) circle (3pt);
\draw [fill=black] (1.875,2.375) circle (3pt);
\draw [color=black] (1.3333333333333333,4) circle (3pt);
\draw [color=black] (3.6,1.8) circle (3pt);
\end{scriptsize}
\end{tikzpicture} 
    \caption{A sub-star $P(4,2)$}
    \label{fig: a,b}
\end{figure}

The next result illustrates how $S(\ell,s,h,c)$ works in practice and enables induction approach. All numbers are assumed to be sufficiently big to enable the specialization step. 
\begin{proposition}[Specialization]\label{prop: specialization}
Let $H$ be a fixed plane in $\PP^3$ and let $Y=B(L,C,(A,B))$ be a scheme in $\PP^3$ consisting of:
\begin{itemize}
    \item $L$ general lines;
    \item $C$ crosses with general vertices in $H$;
    \item a sub-star $P(A,B)$ in $H$.
\end{itemize}
There exists a subscheme $Y'$ of $\PP^3$, which together with $Y$ are members of the same flat family $\mathcal{Y}$ of subschemes of $\PP^3$, such that the trace subscheme $Y'_{tr}$ of $Y$ on $H$ consists of
\begin{itemize}
    \item $\ell+2c+h$ general lines in $H$;
    \item $s$ general double points in $H$;
\end{itemize}
coming from specializations and
\begin{itemize}
    \item $C-c-h$ zero-dimensional subschemes of $H$ of length $2$ supported on $C-c-h$ general points in $H$;
    \item the set $P(A,B)$;
\end{itemize}
already in $H$ and
\begin{itemize}
    \item $L-\ell-2s$ general points.
\end{itemize}
Then $Y_{res}$ is a subscheme in $\PP^3$ consisting of 
\begin{itemize}
    \item $L_{res}=L-\ell+h-2s$ general lines;
    \item $C_{res}=C-c+s-h$ crosses with general vertices in $H$;
    \item a sub-star $P(A_{res},B_{res})$, where $A_{res}=\ell+2c+h$ and $B_{res}=c$.
\end{itemize}
\end{proposition}
\begin{proof}
We perform $S(\ell_{tr},s_{tr},h_{tr},c_{tr})$
and we are done.
\end{proof}
\begin{notation}[Trace and residual schemes]
For the specialization described in Proposition \ref{prop: specialization} and the resulting trace and residual schemes we use the following block notation. For convenience we let extend the notation of the $B$ schemes by including $d$, the degree of considered forms.
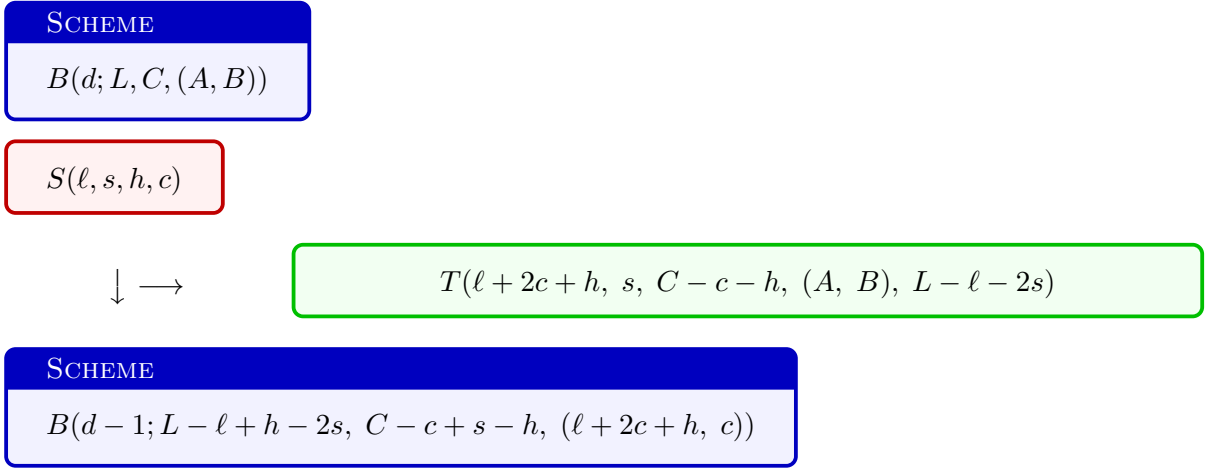
\initialN{$B(d;L,C,(A,B))$}
\reductionN{$S(\ell,s,h,c)$}
{$T(\ell+2c+h,\; s,\; C-c-h,\; (A,\; B),\; L-\ell-2s)$}
{$B(d-1;L-\ell+h-2s,\; C-c+s-h,\; (\ell+2c+h,\; c))$}

Here the numbers in the trace subscheme $T$ represent lines, double points, subschemes of length 2, the next two numbers determine a sub-star, and the last one denotes the number of general points.
\end{notation}

\section{The key reduction lemma and the induction argument}\label{sec: key reduction}
The main idea of the proof of the good postulation of lines in $\PP^3$ is based on the induction argument with respect to the degree of considered forms. We will now show, how to embed Proposition \ref{prop: specialization} effectively in this procedure. 
\begin{lemma}[Reduction Lemma]\label{lem: reduction}
For a positive integer $k$ and $d=3k$ let $Y_{d}$ be a subscheme in $\PP^3$ consisting of 
\begin{itemize}
    \item $\ell=\left\{\begin{array}{ccc}
    \frac{k(k+5)+4}{2}&\mbox{ for }&k\notin 2\Z  \\
    \frac{k(k+5)+2}{2}&\mbox{ for }&k\in 2\Z
    \end{array}\right.$ lines and
    \item $c=\left\{\begin{array}{ccc}
    \frac{k^2-1}{2}&\mbox{ for }&k\notin 2\Z  \\
    \frac{k^2}{2}&\mbox{ for }&k\in 2\Z
    \end{array}\right.$ crosses and
    \item the same number $c$ of points of a sub-star $P(k+1,\frac{k+1}{2})$ points in the odd case and of $P(k+1,\frac{k}{2})$ points in the even case.
\end{itemize}
For $k\geq 4$ there exists a sequence of three specializations $Y_d=X_{0} \to X_{1} \to X_{2} \to X_{3}$ as described in Proposition \ref{prop: specialization} with the residual subschemes $X_{a+1}=(X_{a})_{res}$ for $a=0,1,2$ and such that $X_3=Y_{d-3}$.
\end{lemma}
\begin{proof}
For $k$ odd we perform
\initial{$B\left(3k;\frac{1}{2}k(k+5)+2,\frac{1}{2}(k^2-1),(k+1,\frac{1}{2}(k+1))\right)$}

\reductionL{$S\left(\frac{1}{2}(k-5),0,\frac{1}{2}(k-1),2\right)$}
{$T\left(k+1,0,\frac12 k(k-1)-2,(k+1, \frac12(k+1)),\frac12(k^2+4k+9)\right)$}
{$B\left(3k-1;\frac12 k(k+5)+4,\frac12k(k-1)-2,(k+1,2)\right)$}{red: 3k+2 odd comes here}

\reductionL{$S(0,3,1,\frac{1}{2}(k-1))$}{$T\left(k,3,\frac12(k^2-1)-k-2,(k+1,2),\frac12 k(k+5)-2\right)$}{$B\left(3k-2;\frac12 k(k+5)-1,\frac12(k-1)^2,(k,\frac12 (k-1))\right)$}{red: 3k odd comes here}

\reductionL{$S(0,\frac{1}{2}(k+1),1,\frac{1}{2}(k-1))$}{$T\left(k,\frac12(k+1),\frac12 k(k-3),(k,\frac12(k-1)),\frac12k(k+3)-2\right)$}{$B\left(3(k-1);\frac12(k-1)(k+4)+1,\frac12(k-1)^2,(k,\frac12(k-1))\right)$}{sch: interesting T}

Note that the last displayed scheme is exactly the next displayed scheme with the value of $k$ diminished by $1$. So the case of $k$ odd has been reduced to the case of $k-1$, which is an even number.

In the case $k$ is even we perform
\initial{$B\left(3k;\jd k(k+5)+1,\jd k^2,(k+1,\jd k)\right)$}

\reductionL{$S\left(\jd k-1,0,\jd k,1\right)$}{$T\left(k+1,0,\frac12k(k-1)-1,(k+1,\frac12k),\frac12(k+2)^2\right)$}{$B\left(3k-1;\frac12k(k+5)+2, \frac12k(k-1)-1,(k+1,1)\right)$}{red: 3k+2 even comes here}

\reductionL{$S\left(0,1,0,\jd k\right)$}{$T\left(k,1,\frac12k(k-2)-1,(k+1,1),\frac12k(k+5)\right)$}{$B\left(3k-2;\frac12k(k+5),\frac12k(k-2),(k,\frac12k)\right)$}{red: 3k even comes here}

\reduction{$S\left(0,\jd k,0,\jd k\right)$}{$T\left(k,\frac12k,\frac12k(k-3),(k,\frac12k),\frac12k(k+3)\right)$}{$B\left(3(k-1);\frac12(k-1)(k+4)+2,\frac12((k-1)^2-1),(k,\frac12k)\right)$}

Note again that the last displayed scheme corresponds to the scheme initial in the proof with the value of $k$ diminished by $1$. So the case of $k$ even has been reduced to $k-1$, which is an odd number.

The initial assumption $k\geq 4$ is required in the first reduction step of $k$ odd, which under this assumption is $k\geq 5$. So, the first case not covered by this Lemma is 
\initial{$B(9;14,4,(4,2))$}
for $k=3$.


The correctness of the invariants appearing in each sequence of specializations can be checked by dull calculations or by using the software package we prepared in Singular \cite{Singular} and which is available online at \cite{SingularFiles}.
\end{proof}

\begin{corollary}\label{cor: reduction engine}
Lemma \ref{lem: reduction} shows that starting with arbitrary $k$ and performing $6$ reduction steps, we arrive to the case of $k-2$, which is of the same parity.\\
All linear systems appearing in the Reduction Lemma are expected to be empty. This property would follow for all of them (because of the reduction applied) provided we show that there is no surface of degree $9$ containing the scheme $B(9;14,4,(4,2))$. We deal with this case in Proposition \ref{prop: B(14,4,4,2)}.
\end{corollary}
Now we will show how the subschemes $Y(\ell,p)$ from the Main Theorem can be specialized to subschemes studied in Lemma \ref{lem: reduction} and its proof.
\begin{proposition}[Initial specialization]\label{prop: initial specialization}
All schemes $Y(\ell,p)$ can be specialized to those appearing in Lemma \ref{lem: reduction} in at most $3$ steps.
\end{proposition}
\begin{proof}
We consider cases depending on the congruence of $d$ modulo $3$ and the parity of $k=\lfloor\frac{d}{3}\rfloor$.
\textbf{Case $d=3k$.}
We have $Y\left(\frac{1}{2}k(3k+5)+1,0\right)=B\left(3k;\frac{1}{2}k(3k+5)+1,0,(0,0)\right)$ and we perform
\initial{$B\left(3k;\frac{1}{2}k(3k+5)+1,0,(0,0)\right)$}

\reduction{$S\left(k+1,\frac12k(k-1),(0,0)\right)$}{$T\left(k+1,\frac12k(k-1),0,(0,0),\frac12k(k+5)\right)$}{$B\left(3k-1;\frac12k(k+5),\frac12k(k-1),(k+1,0)\right)$}

Further reductions depend on the parity of $k$. If \textbf{$k$ is odd}, we perform

\setboxcounter{schemeM}{11}
\reductionM{$S\left(0,1,1,\frac12(k-1)\right)$}{$T\left(k,1,\frac12(k^2-2k-1),(k+1,0),\frac12k(k+5)-2\right)$}{$B\left(3k-2;\frac{1}{2}k(k+5)-1,\frac{1}{2}(k-1)^2,(k,\frac12(k-1))\right)$}{odd}

which is Scheme \ref{red: 3k odd comes here} in Lemma \ref{lem: reduction}, so we are done here.

In the case \textbf{$k$ is even} we perform

\setboxcounter{scheme}{10}
\initial{$B\left(3k-1;\frac12k(k+5),\frac12k(k-1),(k+1,0)\right)$}

\setboxcounter{schemeM}{11}
\reductionM{$S\left(0,0,0,\frac12k\right)$}{$T\left(k,0,\frac12k(k-2),(k+1,0),\frac12k(k+5)\right)$}{$B\left(3k-2;\frac12k(k+5),\frac12k(k-2),(k,\frac12k)\right)$}{even}

which is Scheme \ref{red: 3k even comes here} in Lemma \ref{lem: reduction}, so we are done here.
\medskip

\setboxcounter{scheme}{12}
\noindent
\textbf{Case} $d=3k+1$. We have $Y\left(\frac{1}{2}k(3k+7)+2,0\right)=B\left(3k+1;\frac{1}{2}k(3k+7)+2,0,(0,0)\right)$ 
and we perform 
\initial{$B\left(3k+1;\frac{1}{2}k(3k+7)+2,0,(0,0)\right)$}

\reduction{$S\left(k+2,\frac12k(k-3),0,0\right)$}{$T\left(k+2,\frac12k(k-3),0,(0,0),\frac12k(k+11)\right)$}{$B\left(3k;\frac12k(k+11),\frac12k(k-3),(k+2,0)\right)$} 

Further reductions depend on the parity of $k$. If \textbf{$k$ is odd}, we perform

\setboxcounter{schemeM}{14}
\reductionM{$S\left(0,1,1,\frac12(k-1)\right)$}{$T\left(k,1,\frac12(k^2-4k-1),(k+2,0),\frac12k(k+11)-2\right)$}{$B\left(3k-1;\frac12k(k+11)-1,\frac12(k^2-4k+1),(k,\frac12(k-1))\right)$}{odd}

\setboxcounter{schemeM}{15}
\reductionM{$S\left(0,\frac12(3k+1),1,\frac12(k-1)\right)$}{$T\left(k,\frac12(3k+1),\frac12k(k-5), (k, \frac12(k-1)),\frac12k(k+5)-2\right)$}{$B\left(3k-2;\frac12k(k+5)-1,\frac12(k-1)^2,(k,\frac12(k-1))\right)$}{odd}

which is Scheme \ref{red: 3k odd comes here} in Lemma \ref{lem: reduction}, so we are done here.

In the case \textbf{$k$ is even} we perform
\setboxcounter{scheme}{13}
\initial{$B\left(3k;\frac12k(k+11),\frac12k(k-3),(k+2,0)\right)$} 

\setboxcounter{schemeM}{14}
\reductionM{$S\left(0,0,0,\frac12k\right)$}{$T\left(k,0,\frac12k(k-4),(k+2,0),\frac12k(k+11)\right)$}{$B\left(3k-1;\frac12k(k+11),\frac12k(k-4),(k,\frac12k)\right)$}{even}

\setboxcounter{schemeM}{15}
\reductionM{$S\left(0,\frac32k,0,\frac12k\right)$}{$T\left(k,\frac32k,\frac12k(k-5), (k, \frac12k),\frac12k(k+5)\right)$}{$B\left(3k-2;\frac12k(k+5),\frac12k(k-2),(k,\frac12k)\right)$}{even}

which is Scheme \ref{red: 3k even comes here} in Lemma \ref{lem: reduction}, so we are done here.
\medskip

\noindent
\textbf{Case} $3k+2$. We consider first the case when \textbf{$k$ is odd} the first specialization denoted by using $S^*$ is the usual specialization $S$ with all parameters indicated extended by specializing the initial $(k+1)$ collinear points to $H$.

\setboxcounter{scheme}{16}
\initial{$Y(\frac{1}{2}k(3k+9)+3,k+1)$}

\reduction{$S^\mathbf{*}\left(k+2,\frac12k(k-3)-1,0,0\right)$}{$T\left(k+2,\frac12k(k-3)-1,0,(0,0),\frac12k(k+13)+3+\mathbf{(k+1)}\right)$}{$B\left(3k+1;\frac12k(k+13)+3,\frac12k(k-3)-1,(k+2,0)\right)$}

\reduction{$S\left(0,k+2,1,\frac12(k-1)\right)$}{$T\left(k,k+2,\frac12(k^2-4k-3),(k+2,0),\frac12k(k+9)-1\right)$}{$B\left(3k;\frac12k(k+9),\frac12(k^2+1)-k,(k,\frac12(k-1))\right)$}

\reduction{$S\left(k-3,\frac12(k-1),0,2\right)$}{$T\left(k+1,\frac12(k-1),\frac12(k^2+1)-k-2,(k,\frac12(k-1)),\frac12k(k+5)+4\right)$}{$B\left(3k-1;\frac12k(k+5)+4,\frac12k(k-1)-2,(k+1,2)\right)$}
This is Scheme \ref{red: 3k+2 odd comes here}  in Lemma \ref{lem: reduction}, so we are done.
\medskip

Finally, for \textbf{$k$ even} in case we have first the $S^*$ specialization. 

\initial{$Y(\frac{1}{2}k(3k+9)+3,k+1)$}

\reduction{$S^\mathbf{*}\left(k+2,\frac12k(k-3)-1,0,0\right)$}{$T\left(k+2,\frac12k(k-3)-1,0,(0,0),\frac12k(k+13)+3+\mathbf{(k+1)}\right)$}{$B\left(3k+1;\frac12k(k+13)+3,\frac12k(k-3)-1,(k+2,0)\right)$}

\reduction{$S\left(0,k+1,0,\frac12k\right)$}{$T\left(k,k+1,\frac12k(k-4)-1,(k+2,0),\frac12k(k+9)+1\right)$}{$B\left(3k;\frac12k(k+9)+1,\frac12k(k-2),(k,\frac12k)\right)$}

\reduction{$S\left(k-1,\frac12k,0,1\right)$}{$T\left(k+1,\frac12k,\frac12k(k-2)-1,(k,\frac12k),\frac12k(k+5)+2\right)$}{$B\left(3k-1;\frac12k(k+5)+2,\frac12k(k-1)-1,(k+1,1)\right)$}
which is Scheme \ref{red: 3k+2 even comes here} considered in Lemma \ref{lem: reduction}.

This ends the whole proof.
\end{proof}

\begin{corollary}\label{cor: number of lines 12}
Proposition \ref{prop: initial specialization} implies that all schemes $Y(\ell,p)$ from the Main Theorem \ref{thm: main}, under the condition $d\geq 12$ can be reduced to schemes of the type $B(L,C,(A,B))$.    
\end{corollary}

For the remaining cases see the next Section.

\section{The cases not covered by the induction}\label{sec: preinduction cases}
We begin with the case crucial for the reduction engine described in Lemma \ref{lem: reduction}.
\begin{proposition}\label{prop: B(14,4,4,2)}
The linear system of surfaces of degree $9$ containing the scheme $B(9;14,4,(4,2))$ is empty.    
\end{proposition}
\begin{proof}
In the proof we run a reduction process, which involves the same type of specializations as the Reduction Lemma \ref{lem: reduction} but for numerical reason we need to apply them in a more subtle way. We begin with
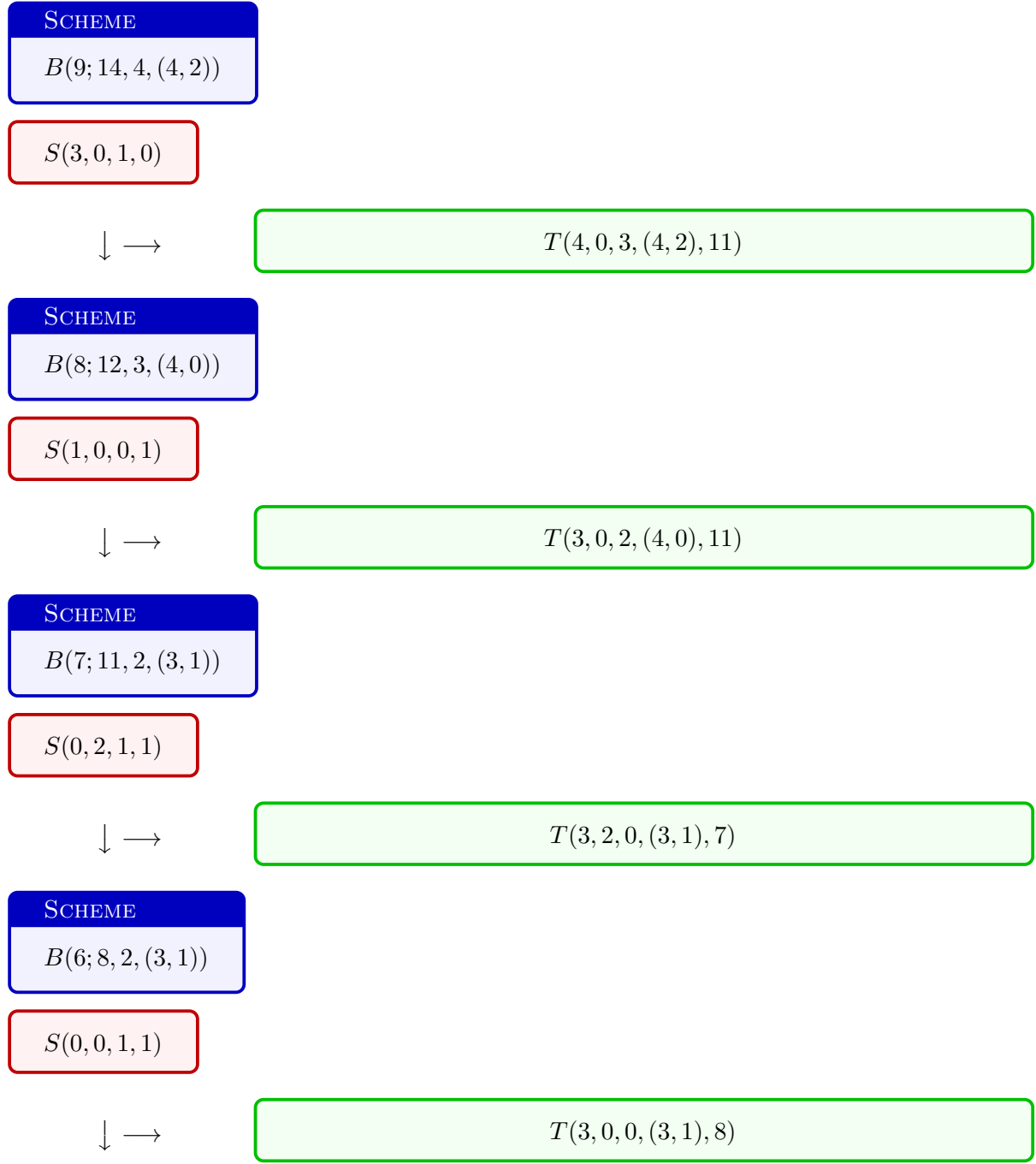
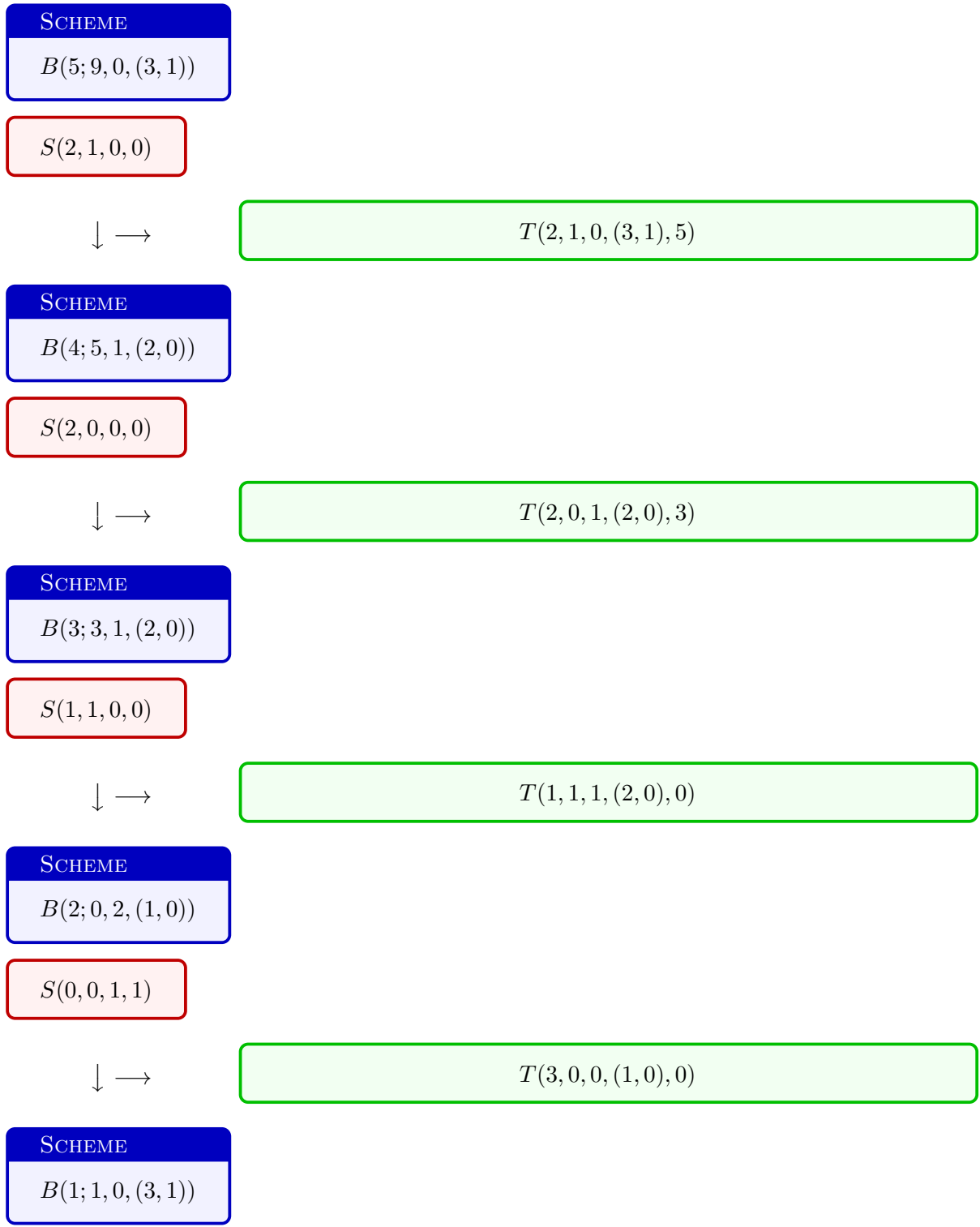
\initialN{$B(9;14,4,(4,2))$}
\reductionN{$S(3,0,1,0)$}{$T(4,0,3,(4,2),11)$}{$B(8;12,3,(4,0))$}
\reductionN{$S(1,0,0,1)$}{$T(3,0,2,(4,0),11)$}{$B(7;11,2,(3,1))$}
\reductionN{$S(0,2,1,1)$}{$T(3,2,0,(3,1),7)$}{$B(6;8,2,(3,1))$}
\reductionN{$S(0,0,1,1)$}{$T(3,0,0,(3,1),8)$}{$B(5;9,0,(3,1))$}
\reductionN{$S(2,1,0,0)$}{$T(2,1,0,(3,1),5)$}{$B(4;5,1,(2,0))$}
\reductionN{$S(2,0,0,0)$}{$T(2,0,1,(2,0),3)$}{$B(3;3,1,(2,0))$}
\reductionN{$S(1,1,0,0)$}{$T(1,1,1,(2,0),0)$}{$B(2;0,2,(1,0))$}
\reductionN{$S(0,0,1,1)$}{$T(3,0,0,(1,0),0)$}{$B(1;1,0,(3,1))$}
Note that the sub-star $(3,1)$ is just $2$ general points in the plane $H$. As there is obviously no plane vanishing along a line and two general points, we are done.
\end{proof}

\begin{remark}
By Corollary \ref{cor: number of lines 12} our proof of the Main Theorem works for forms of degree at least $12$. While the arguments can be ad hoc changed to include cases of lower values of $d$ we content ourselves with referring to the original work of Hartshorne and Hirschowitz \cite{HarHir82}. The main purpose of our note is to illustrate a new method.
\end{remark}

\section{Good postulation in $\PP^2$}\label{sec: postulation P2}
In the course of reductions outlined in sections \ref{sec: main tools} and \ref{sec: key reduction} we need to handle postulation of trace schemes on forms in $\PP^2$ (the schemes denoted by $T$ and appearing the green boxes). 

Our idea here follows the pattern developed for subschemes in $\PP^3$, namely we specialize elements of subschemes $T$ to a line $L$ in $\PP^2$. We have the following specializations:
\begin{itemize}
    \item double point (the residual scheme is a simple point);
    \item the subscheme of length $2$ embedded in the line (the residual scheme is empty);
    \item points from a line in a sub-star (in this case we specialize everything to this particular line);
    \item a general point.
\end{itemize}
Of course, on $L$ all zero-dimensional subschemes are divisors, so the good postulation is automatic.

It convenient to introduce the following notation.
\begin{notation}
We will denote by $V(d,h,(a,b),p)$ a subscheme of $\PP^2$ consisting of
\begin{itemize}
    \item $d$ double points;
    \item $h$ subschemes of lengths $2$;
    \item an $(a,b)$ sub-star and
    \item $p$ general points.
\end{itemize}
\end{notation}
\begin{theorem}\label{thm: good postulation in P2}
All trace schemes appearing in the reduction steps in Sections \ref{sec: main tools} and \ref{sec: key reduction} have good postulation.
\end{theorem}
\begin{proof}
As handling all possibilities would be very long, we content ourselves with explaining two examples. All remaining cases are left to the motivated reader.

We begin with the trace scheme
$$T\left(k+1,0,\frac12 k(k-1)-2,(k+1, \frac12(k+1)),\frac12(k^2+4k+9)\right)$$
appearing right before Scheme \ref{red: 3k+2 odd comes here}.

We claim in effect that there is no form of degree $2k-1$ vanishing on 
$$V\left(0,\jd k(k-1)-2,(k+1,\jd(k+1)),\jd(k^2+4k+9)\right).$$

In the first step we specialize to one of sub-star lines. The line contains $k-1$ points and we specialize additional $k+1$ general points on that line. As there is no form of degree $2k-1$ on this line vanishing at all $2k=(k-1)+(k+1)$ points, the line must be a component of any form of degree $2k-1$ in $\PP^2$ vanishing along the given scheme, so we can pass to looking at the residual scheme, which is
$$V\left(0,\jd k(k-1)-2,(k,\jd(k-1)),\jd(k^2+2k+7)\right).$$

and we claim now that no form of degree $2k-2$ vanishes on that scheme.

For the next step we note that in the sub-star $(k,\jd(k-1))$ there is just one line still containing $k-1$ points. We specialize to this line additional $k$ general points, so that there are $2k-1$ points altogether. Again this line must be a component of any form of degree $2k-2$ on $\PP^2$ vanishing on $V$. So we need to study the residual subscheme, which is 
$$V\left(0,\jd k(k-1)-2,(k-1,\jd(k-1)),\jd(k^2+7)\right).$$
Now we will reduce the number of subschemes of length $2$. First we specialize $(k-1)$ of them on a single line and so that their restriction to that line has length $2(k-1)$. Again this line must then a component of any curve of degree $2k-3$ in $\PP^2$ vanishing on $V$. The residual scheme is now
$$V\left(0,\jd k(k-3)-1,(k-1,\jd(k-1)),\jd(k^2+7)\right).$$
In the last step we specialize $(k-2)$ subschemes of lengts $2$ and one general point. By the same token as seen above we arrive to study forms of degree $2k-5=2(k-2)-1$ vanishing along
$$V\left(0,\jd k(k-5)+1,(k-1,\jd(k-1)),\jd(k^2+5)\right).$$
But this is exactly the scheme we began with with parameter $(k-2)$ in place of $k$.

It remains to show that the reductions lead to a subscheme with good postulation. 
As we started our argument with $k$ odd and we are able to lower its value by $2$ applying the above steps over and over, we will finally come to the value $k=3$. Here we have $V(0,1,(4,2),15)$. We need to show that there is no form of degree $5$ vanishing along this subscheme. Specializing to one of the sub-star lines
\begin{itemize}
    \item $1$ subscheme of length $2$ and
    \item $2$ general points,
\end{itemize}
we obtain as the residual scheme $V(0,0,(3,1),13)$.
Recalling that a $(3,1)$ sub-star is just $2$ general points, we are reduced to the question of the existence of a form of degree $4$ vanishing in $15$ general points, which is well-known.
\medskip

As the second example we study the trace scheme
$$T\left(k,\frac12(k+1),\frac12 k(k-3),(k,\frac12(k-1)),\frac12k(k+3)-2\right)$$
appearing right before Scheme \ref{sch: interesting T}. It contains double points so something new is needed.

We claim that there is no form of degree $2(k-1)$ vanishing along
$$V=\left(\frac12(k+1),\frac12 k(k-3),(k,\frac12(k-1)),\frac12k(k+3)-2\right).$$
In the first step we take a general line and specialize there one double point $2P$, $(k-2)$ schemes of length $2$ and one general point. This forces the line to be contained in the curve of degree $2k-2$ vanishing along $V$, should it exist.

The residual scheme is
$$V=\left(\frac12(k-1),\frac12 k(k-5)+2,(k,\frac12(k-1)),\frac12k(k+3)-3\right) \mbox{ plus a single point }P.$$
Now, we take a general line through $P$ and put on that line additionally $(k-3)$ subschemes of length $2$ and $3$ general points. Again the conditions on the line are so that we can pass to the residual scheme
$$V=\left(\frac12(k-1),\frac12 k(k-7)+5,(k,\frac12(k-1)),\frac12k(k+3)-6\right).$$
We specialize now to the only line with $(k-1)$ points in the sub-star $k-2$ general points. This forces the line to be a component of a curve of degree $2(k-2)$ vanishing along $V$.

The residual scheme is
$$V=\left(\frac12(k-1),\frac12 k(k-7)+5,(k-1,\frac12(k-1)),\frac12k(k+1)-4\right).$$
Specializing now to a line in the sub-star $k-1$ general points, we arrive finally to the subscheme
$$V=\left(\frac12(k-1),\frac12 k(k-7)+5,(k-2,\frac12(k-3)),\frac12k(k-1)-3\right).$$

It easy to check that this scheme is exactly the scheme we started with with $k$ replaced by $k-2$.

Similarly as before we need to consider the case $k=3$ for which we get $V(2,0,(3,1),7)$. Since the sub-star $(3,1)$ is just two general points, we need to show that there is no plane quartic passing through $2$ general double points and $9$ general points, which is well known.

\end{proof}

\paragraph{Acknowledgements.} The present work was initiated during a Workshop on Algebraic Geometry held in Lanckorona in October 2024. We thank Willa Tadeusz for hospitality and the Excellence Small Working Groups program DNa.711/IDUB/ESWG/2024/01/00010 for financial support.

\bigskip
\noindent
Marcin Dumnicki and Halszka Tutaj-Gasińska,\\
Department of Mathematics,
Jagiellonian University,
Łojasiewicza 6,
PL-30-348 Kraków, Poland\\
\nopagebreak
\textit{E-mail address:} \texttt{marcin.dumnicki@uj.edu.pl}\\
\textit{E-mail address:} \texttt{halszka.tutaj-gasinska@uj.edu.pl}\\

\noindent
Mikołaj Le Van, Grzegorz Malara, Tomasz Szemberg and Justyna Szpond,\\
Department of Mathematics,
University of the National Education Commission Krakow,
Podchor\c a\.zych 2,
PL-30-084 Krak\'ow, Poland. \\
\nopagebreak
\textit{E-mail address:} \texttt{mikolaj.levan@student.uken.krakow.pl}\\
\textit{E-mail address:} \texttt{grzegorz.malara@gmail.com}\\
\textit{E-mail address:} \texttt{tomasz.szemberg@gmail.com}\\
\textit{E-mail address:} \texttt{szpond@gmail.com}\\
\end{document}